\documentclass[psamsfonts]{amsart}

\usepackage{amssymb,amsfonts}
\usepackage[all,arc]{xy}
\usepackage{enumerate}
\usepackage{mathrsfs}
\usepackage{graphicx}

\newtheorem{theorem}{Theorem}[section]
\newtheorem{corollary}[theorem]{Corollary}
\newtheorem{proposition}[theorem]{Proposition}

\theoremstyle{definition}
\newtheorem{definition}[theorem]{Definition}

\newtheorem{example}[theorem]{Example}

\newtheorem{remark}[theorem]{Remark}

\let\phi=\varphi

\def\Ann{\operatorname{Ann}}

\def\Min{\operatorname{Min}}
\def\Spec{\operatorname{Spec}}

\let\oldbigwedge\bigwedge
\def\BIGwedge{{\textstyle\oldbigwedge}}
\def\medwedge{{\scriptstyle\oldbigwedge}}
\def\bigwedge{\mathchoice{\BIGwedge}{\BIGwedge}{\medwedge}{}}

\DeclareMathOperator{\Nil}{Nil}

\DeclareMathOperator{\Id}{Id}

\DeclareMathOperator{\Stone}{Stone}
\DeclareMathOperator{\pcomp}{pcomp}
\DeclareMathOperator{\Skel}{Skel}
\DeclareMathOperator{\Dns}{Dns}

\let\hat=\widehat

\let\epsilon=\varepsilon

\begin{document}
\title{Pseudocomplementation and Minimal Prime Ideals in Semirings}

\author{Peyman Nasehpour}
\address{Peyman Nasehpour, Department of Engineering Science, Faculty of Engineering, University of Tehran, Tehran, Iran}

\email{nasehpour@gmail.com}



\thanks{Dedicated to Professor Carlos Mart\'{\i}n-Vide}

\begin{abstract} In the first section of the present work, we introduce the concept of pseudocomplementation for semirings and show semiring version of some known results in lattice theory. We also introduce semirings with pc-functions and prove some interesting results for minimal prime ideals of such semirings. In the second section, some classical results for minimal prime ideals in ring theory are generalized in the context of semiring theory.

\end{abstract}

\maketitle

\section{Introduction}

A semirings is a ring-like structure, where subtraction is either impossible or disallowed. Commutative semirings with nonzero identity are important ring-like structures with so many applications in science and engineering (\cite[p. 225]{Golan2003}) and are considered to be interesting generalizations of bounded distributive lattices and commutative rings with nonzero identities (\cite[Example 1.5]{Golan1999}). The concept of complemented elements in semirings have been investigated in chapter five of the book \cite{Golan1999}. As a matter of fact, complemented elements play an important part in the semiring representation of the semantics of computer programs (\cite{ManesArbib1986}). On the other hand, the concept of pseudocomplementation is a well-developed notion in lattice theory (\cite{RasiowaSikorski1963}). Since pseudocomplementation has been recently defined and developed for other algebraic structures like semigroups with zero (\cite{JacksonStockes2004} and \cite{Cirulis2007}) and has important applications in computer science (\cite{DuntschWinter2006} and \cite{Guttmann2016}), it seems quite natural and useful if the concept of pseudocomplementation can be defined, developed, and investigated in the context of semiring theory, the task that we will try to do in \S 1.

Note that since different authors use the term ``semiring'' for different concepts, it is crucial to clarify what it is meant by semiring in this work. More on semirings can be found in the books \cite{Golan1999} and \cite{HebischWeinert1998}.

In this work, by a semiring, we understand an algebraic structure, consisting of a nonempty set $S$ with two operations of addition and multiplication such that the following conditions are satisfied:

\begin{enumerate}
\item $(S,+)$ is a commutative monoid with identity element $0$;
\item $(S,\cdot)$ is a commutative monoid with identity element $1 \not= 0$;
\item Multiplication distributes over addition, i.e. $a\cdot (b+c) = a \cdot b + a \cdot c$ for all $a,b,c \in S$;
\item The element $0$ is the absorbing element of the multiplication, i.e. $s \cdot 0=0$ for all $s\in S$.
\end{enumerate}

Interesting examples for semirings include the tropical algebra $(\mathbb T, \max, +)$, with tropical numbers $\mathbb T = [-\infty, +\infty)$, which is fundamental for the calculations in tropical geometry (\cite{ItenbergMikhalkinShustin2009} and \cite{MikhalkinRau2015}) and $(\hat{\mathbb R}, \min, +)$, with $\hat{\mathbb R} = (-\infty, +\infty]$, which has applications in shortest path problems (\cite{GondranMinoux2008}).

By an ordered semiring, we mean a semiring $(S, +, \cdot)$ with a partial order $\leq$ on $S$ such that the following conditions hold:

\begin{enumerate}

\item If $s \leq t$, then $s+u \leq t+u$ for any $s,t,u \in S$;

\item If $s \leq t$ and $0 \leq u$, then $su \leq tu$ for any $s,t,u \in S$.

\end{enumerate}

An ordered semiring is called to be positive, if $0$ is its least element, i.e. $0 \leq s$ for all $s\in S$. For example $(\Id(S),+,\cdot, \subseteq)$ is a positive semiring, where $S$ is itself an arbitrary semiring (see Proposition \ref{positivepcsemiring}). For more on ordered semirings, one may refer to chapter 2 of the book \cite{Golan2003}.

In \S 1, we introduce the pseudocomplemented, stone, and dense elements in semirings and prove some nice results related to these elements, similar to what we have in lattice theory (See Proposition \ref{pseudo1}, Proposition \ref{stone1}, Proposition \ref{simple}, Proposition \ref{stone2}, Theorem \ref{skel1}, and Proposition \ref{dense1}).

Let us recall that a nonempty subset $I$ of a semiring $S$ is called an ideal, if $a+b \in I$ and $sa \in I$ for all $a,b \in I$ and $s\in S$. An ideal $P \neq S$ is defined to be a prime ideal of $S$, if $ab\in P$ implies either $a\in P$ or $b\in P$. A prime ideal $P$ of a semiring $S$ is called to be a minimal prime ideal of $S$, if $I \subseteq P$ implies either $I = (0)$ or $I = P$. For more on ideals of a semiring, one can refer to chapter 6 and chapter 7 of the book \cite{Golan1999}.

In fact, in \S 1, we do more and introduce pc-functions (See Definition \ref{pc-function}) in this way that in a semiring $S$, we define a function $*: S \longrightarrow S$ to be a pseudocomplemented function (for short pc-function), if $s \cdot *(s) = 0$ for all $s\in S$ - for the ease of calculation, we denote $*(s)$ by $s^*$ - and for example in Theorem \ref{pcfunctionminimalprime2}, we prove that if $S$ is a semiring and $*: S \longrightarrow S$ a pc-function such that $0^*=1$ and $(s+s^*)^* = 0$ for any $s\in S$, and $P$ is a prime ideal of $S$, then the following statements are equivalent:

\begin{enumerate}

\item If $s\in P$, then $s^* \notin P$, for each $s\in S$,

\item If $s\in P$, then $s^*{^*} \in P$, for each $s\in S$,

\item $P \cap \{s\in S: s^* = 0 \} = \emptyset$.

\end{enumerate}

Moreover, if one of the above conditions hold, $P$ is a minimal prime ideal of $S$. Those, who are familiar with lattice theory, are aware of this point that these are a generalization of some interesting results for minimal primes in lattice theory (See Theorem \ref{pcfunctionminimalprime}, Corollary \ref{minimalprimepseudo}, and Theorem \ref{pcfunctionminimalprime2}).

Let us recall that if $R$ is a commutative ring with a nonzero identity, then a prime ideal $P$ is a minimal prime ideal of an ideal $I$ in $R$ if and only if for each $x\in P$, there is a $y\notin P$ and a nonnegative integer $i$ such that $yx^i \in I$ (Check Theorem 2.1 in \cite{Huckaba1988}). This classical result in commutative algebra has some interesting corollaries for reduced rings (See corollaries of Theorem 2.1 in \cite{Huckaba1988}). In \S 2, we prove the semiring version of this classical result and its corollaries for nilpotent-free semirings, i.e., semirings with no non-trivial multiplicatively nilpotent elements (Check Theorem \ref{minimalprimehuckaba}, Corollary \ref{minimalprimehuckaba2}, and Corollary \ref{minimalprimehuckaba3}). We end this work by characterizing minimal primes of pseudocomplemented semirings (See Theorem \ref{minimalpbdl}).

\section{Pseudocomplemented, Stone and Dense Elements in Ordered Semirings}

Let us recall that in a meet-semilattice $L$ with zero, an element $a^*$ is a pseudocomplement of $a$ ($\in L$), if $a \wedge a^* = 0$ and $a \wedge x = 0$ implies that $x \leq a^*$ for any $x\in L$ (\cite[Sect. 6.2]{Gratzer2011}). The pseudocomplement elements of semigroups with zero are defined in \cite{JacksonStockes2004} similarly. We define the concept of pseudocomplemented elements for ordered semirings as follows:

\begin{definition}

\label{pseudodef}

Let $S$ be an ordered semiring. We say an element $s\in S$ is a pseudocomplemented element of $S$ if there exists an element $s_p\in S$ such that the following properties hold:

\begin{enumerate}

\item $ss_p = 0$,

\item If $sx = 0$, then $x \leq s_p$ for each $x\in S$.

\end{enumerate}

If such a $s_p$ exists, then by definition, we call it the pseudocomplement of $s$.

\end{definition}

Note that any element in an ordered semiring $S$ has at most one pseudocomplement and since the only element that annihilates $1$ is $0$, the element $1$ is a pseudocomplemented element of $S$ and its pseudocomplement is $0$.

Let us recall that the skeleton of a meet-semilattice $L$ is defined to be the set $\Skel(L) = \{s^* : s\in L \}$ (\cite[Sect. 6.2]{Gratzer2011}). We similarly define the skeleton of an ordered semiring as follows:

\begin{definition}
Let $S$ be an ordered semiring. If an element $s\in S$ possesses a pseudocomplement, we denote its pseudocomplement by $s^*$. We collect all the pseudocomplemented elements of $S$ in the set $\pcomp(S)$. We define a semiring $S$ to be pseudocomplemented if $S = \pcomp(S)$. We define the skeleton of the semiring $S$ to be the set $\Skel(S) = \{s^* : s\in \pcomp(S) \}$.
\end{definition}

\begin{proposition}

\label{positivepcsemiring}

If $S$ is a semiring, then $(Id(S),+,\cdot, \subseteq)$ is a positive and pseudocomplemented semiring and the pseudocomplement of the ideal $I$ is the ideal $K$ generated by all ideals $J$ such that $I\cdot J = (0)$.

\begin{proof}
This point that $\Id(S)$ is a positive semiring is straightforward. We only show that each element of $\Id(S)$ is pseudocomplemented. For doing so, assume that $I$ is an ideal of $S$, then the set of all ideals of $S$ that annihilates $I$ is nonempty, since $I$ can be annihilated by the zero ideal.

On the other hand, if we set $\Delta =\{J\in \Id(S): I\cdot J = (0)\}$, then the ideal $K$, generated by all ideals $J \in \Delta$, annihilates $I$ and the reason is that any element of $K$ is annihilated by any element of $I$. Also note that $K$ contains all the elements of $\Delta$ and therefore $K$ is the pseudocomplement of $I$, i.e., $K = I^*$ and the proof is complete.
\end{proof}

\end{proposition}

Note that every bounded distributive lattice is a (commutative) semiring. Now we give the following example:

\begin{example}[Example of a positive and pseudocomplemented semiring that is not a pseudocomplemented bounded distributive lattice]

Let $R$ be a commutative ring with a nonzero identity and $\mathfrak{a}, \mathfrak{b}$ and $\mathfrak{c}$ be ideals of $R$ such that $\mathfrak{a} \not= \mathfrak{a}^2$, $\mathfrak{a} \supseteq \mathfrak{b}$, $\mathfrak{a} \supseteq \mathfrak{c}$ and $\mathfrak{b}\cdot \mathfrak{c} = 0$. Then it is obvious that $\mathfrak{a} + \mathfrak{b}\cdot \mathfrak{c} = \mathfrak{a}$, while $(\mathfrak{a}+\mathfrak{b})\cdot (\mathfrak{a}+\mathfrak{c}) = \mathfrak{a}^2$ and in this case $+$ is not distributed over $\cdot$, which means that $(Id(S),+,\cdot, \subseteq)$ is a positive and pseudocomplemented semiring and also a bounded, but not a distributive lattice. As as example, set $S=\Id(\mathbb Z_{n^3})$ to be the semiring of all ideals of the ring $\mathbb Z_{n^3}$, where $n\geq 2$ is a natural number. It is clear that by Proposition \ref{positivepcsemiring}, $S$ is a positive and pseudocomplemented semiring. But if we set $\mathfrak{a} = (n)$ and $\mathfrak{b} = \mathfrak{c} = (n^2)$, then $$\mathfrak{a} + \mathfrak{b}\cdot \mathfrak{c} \ne (\mathfrak{a}+\mathfrak{b})\cdot (\mathfrak{a}+\mathfrak{c}).$$

\end{example}

\begin{proposition}

\label{pseudo1}

Let $S$ be a positive semiring. Then the following statements hold:

\begin{enumerate}

\item If $0 \in \pcomp(S)$, then $0\leq s \leq 0^*$ for any $s\in S$.

\item If $s \in \pcomp(S)$, then $ss^* = 0$. In addition, if $s^* \in \pcomp(S)$, then $s^* s^*{^*} = 0$.

\item If $t\in \pcomp(S)$, then $st = 0$ if and only if $s \leq t^*$.

\item If $s,s^* \in \pcomp(S)$, then $s \leq s^*{^*}$.

\item If $s,s^*, s^*{^*} \in \pcomp(S)$, then $s^*{^*}{^*} = s^*$.

\item If $s, s^* \in \pcomp(S)$, then $st = 0$ if and only if $s^*{^*} t= 0$.

\item If $s,t \in \pcomp(S)$, then $s\leq t$ implies that $t^* \leq s^*$.

\item If $s,s^* \in \pcomp(S)$, then $s\in \Skel(S)$ if and only if $s^*{^*} = s$.

\end{enumerate}

\begin{proof}

(1): Since $s\cdot 0 = 0$ for any $s \in S$, we have $s \leq 0^*$.

(2): It is just a result of the definition of pseudocomplemented elements in a semiring.

(3): Let $t\in \pcomp(S)$. It is clear that if $st = 0$, then $s \leq t^*$. On the other hand, if $s \leq t^*$, then $st \leq tt^* = 0$. But $0$ is the least element of $S$, so $st = 0$.

(4): $s^*{^*}$ is the largest element of $S$ that annihilates $s^*$. But $s^* s = 0$, so $s \leq s^*{^*}$.

(5): By (4), $s^* \leq s^*{^*}{^*}$. Also since $s \leq s^*{^*}$, by (4), we have that $s s^*{^*}{^*} \leq s^*{^*} s^*{^*}{^*} = 0$. Therefore $s^*{^*}{^*} \leq s^*$.

(6): If $st = 0$, then $t \leq s^*$. Now $ s^*{^*} t \leq s^*{^*} s^* = 0$. Conversely, if $s^*{^*} t= 0$, then $st \leq s^*{^*} t= 0$.

(7): If $s \leq t$, then $s t^* \leq t t^* = 0$. This means that $t^* \leq s^*$.

(8): Let $s\in \Skel(S)$. So there is a $t\in \pcomp(S)$ such that $s= t^*$. Since $s,s^*\in \pcomp(S)$, we have that $t^*, t^*{^*} \in \pcomp(S)$. Therefore by (5), we have that $s^*{^*} = t^*{^*}{^*} = t^* = s$. Conversely, let $s^*{^*} = s$. Since $s^* \in \pcomp(S)$ and $s=(s^*)^*$, $s\in \Skel(S)$.
\end{proof}

\end{proposition}

This property that any pseudocomplemented element $s$ of an ordered semiring $S$, is annihilated by its pseudocomplement $s^*$, i.e., $ss^* = 0$, is very interesting and can be considered in a more general context, as we will see in Theorem \ref{pcfunctionminimalprime}, Corollary \ref{minimalprimepseudo}, Theorem \ref{pcfunctionminimalprime2}, and also in Theorem \ref{minimalpbdl}. Based on this property, we generalize this concept and define pseudocomplemented functions (for short pc-functions) as follows:

\begin{definition}

\label{pc-function}

Let $S$ be a semiring. We define a function $*: S \longrightarrow S$ to be a pseudocomplemented function (for short pc-function), if $s \cdot *(s) = 0$ for all $s\in S$. For the ease of calculation, we denote $*(s)$ by $s^*$.

\end{definition}

\begin{theorem}

\label{pcfunctionminimalprime}

Let $S$ be a semiring and $*: S \longrightarrow S$ a pc-function. If a prime ideal $P$ of $S$ has this property that $s\in P$ implies that $s^* \notin P$, then the following statements hold:

\begin{enumerate}

\item If $s\in P$, then $s^*{^*} \in P$.

\item $P$ is a minimal prime ideal of $S$.

\end{enumerate}

\begin{proof}

(1): Let $s\in P$. By assumption, $s^* \notin P$. But $s^* s^*{^*} = 0 \in P$. Since $P$ is prime, $s^*{^*} \in P$.

(2): Suppose $Q \subseteq P$, where $Q$ is a prime ideal of $S$. If there is some $s\in P-Q$, then $s^* \notin P$ and so $s^* \notin Q$. But $ss^* = 0\in Q$, which implies that $s\in Q$ by primeness of $Q$, a contradiction. So $Q = P$ and $P$ is a minimal prime ideal of $S$.
\end{proof}

\end{theorem}

\begin{corollary}

\label{minimalprimepseudo}

Let $S$ be a pseudocomplemented semiring. If a prime ideal $P$ of $S$ has this property that $s\in P$ implies that $s^* \notin P$, then the following statements hold:

\begin{enumerate}

\item If $s\in P$, then $s^*{^*} \in P$.

\item $P$ is a minimal prime ideal of $S$.

\end{enumerate}

\end{corollary}

\begin{theorem}

\label{pcfunctionminimalprime2}

Let $S$ be a semiring and $*: S \longrightarrow S$ a pc-function such that $0^*=1$ and $(s+s^*)^* = 0$ for any $s\in S$. If $P$ is a prime ideal of $S$, then the following statements are equivalent:

\begin{enumerate}

\item If $s\in P$, then $s^* \notin P$, for each $s\in S$,

\item If $s\in P$, then $s^*{^*} \in P$, for each $s\in S$,

\item $P \cap \{s\in S: s^* = 0 \} = \emptyset$.

\end{enumerate}

\begin{proof}

$(1) \Rightarrow (2)$: Theorem \ref{pcfunctionminimalprime}.

$(2) \Rightarrow (3)$: Let $s\in P \cap \{s\in S: s^* = 0 \}$, for some $s\in S$. So $s^*=0$ and therefore $s^*{^*} = 0^* = 1$, which implies that $1\in P$, a contradiction.

$(3) \Rightarrow (1)$: Let $s, s^* \in P$ for some $s\in S$. So $s+s^* \in P$. By assumption, $ (s+s^*)^* = 0$, which means that $P \cap \{s\in S: s^* = 0 \} \neq \emptyset$.
\end{proof}

\end{theorem}

Let us recall that a pseudocomplemented lattice $L$ is called a Stone lattice, if it satisfies the Stone identity: $a^* \vee  a^*{^*} = 1$ for any $a\in L$(\cite{Gratzer2011}). We are inspired to define stone elements in semirings similarly:

\begin{definition}

\label{stonedef}

Let $S$ be a positive semiring.

\begin{enumerate}

\item We define $s\in S$ to be a Stone element of $S$ if $s, s^* \in \pcomp(S)$ and $s^* +  s^*{^*} = 1$. We denote the set of all stone elements of a semiring $S$ by $\Stone(S)$.

\item We define a semiring $S$ to be a Stone semiring, if $S = \Stone(S)$.

\end{enumerate}

\end{definition}

\begin{proposition}

\label{stone1}

Let $S$ be a positive semiring. Then the following statements hold:

\begin{enumerate}

\item If $s\in \Stone(S)$, then $s^*$ is multiplicatively idempotent.

\item If $s\in \Stone(S)$, then $s^2 \leq s$.

\item If $s \in \Stone(S) \cap \Skel(S)$, then $s$ is multiplicatively idempotent.

\end{enumerate}

\begin{proof}

(1): Since $s\in \Stone(S)$, we have that $s^* = s^* \cdot 1 = s^* (s^* +  s^*{^*}) = s^* s^* + s^* s^*{^*} = s^* s^*$.

(2): By Proposition \ref{pseudo1}, we know that $s \leq s^*{^*}$. Therefore $s^2 \leq s s^*{^*}$. But since $s$ is a Stone element of $S$, we have that $ s = s \cdot 1 = s (s^* +  s^*{^*}) = s s^* + s s^*{^*} = s s^*{^*}$.

(3): Since $s \in \Stone(S) \cap \Skel(S)$, we have that $s = s \cdot 1 = s (s^* +  s^*{^*}) = s s^* + s s^*{^*} = ss$.
\end{proof}

\end{proposition}

One of the most simple questions that one may ask about the Stone elements of a semiring $S$ is that when $1\in \Stone(S)$. Surprisingly, this is equivalent to the semiring $S$ to be a simple semiring, i.e., a semiring that for all its elements $s$, $1+s=1$, as we show in the following:

\begin{proposition}

\label{simple}

Let $S$ be a positive semiring and $0 \in \pcomp(S)$. Then the following statements are equivalent:

\begin{enumerate}

\item $1\in \Stone(S)$,

\item $0^* = 1$,

\item $1$ is the largest element of $S$,

\item $S$ is a simple semiring.

\end{enumerate}

\begin{proof}
$(1) \Leftrightarrow (2)$: Let $0^* = 1$. This means that $1\in \Skel(S)$ and so $1= 1^*{^*} = 1^*{^*} + 0 = 1^*{^*} + 1^*$. Obviously $ 1 \cdot 1^* = 0$. Therefore $1\in \Stone(S)$. Conversely, let $1\in \Stone(S)$. So $1^* + 1^*{^*} = 1$ and since $1^* = 0$, we have that $1^*{^*} = 1$. Finally $0^* = 1^*{^*} = 1$.

$(2) \Leftrightarrow (3)$: Straightforward by Proposition \ref{pseudo1}.

$(3) \Leftrightarrow (4)$: Since $0$ is the least element of $S$, if $1$ is the largest element of $S$, we have that $1 \leq s+1 \leq 1$ for any $s\in S$. This means that $S$ is a simple semiring. Conversely, if $S$ is a simple semiring, then $s \leq s+1 = 1$, since $0$ is the least element. This means that $1$ is the largest element of $S$ and the proof is complete.
\end{proof}

\end{proposition}

\begin{proposition}

\label{bdl}

Let $S$ be a multiplicatively idempotent and a positive semiring such that $0 \in \pcomp(S)$. Then the following statements are equivalent:

\begin{enumerate}

\item $1\in \Stone(S)$,

\item $0^* = 1$,

\item $1$ is the largest element of $S$,

\item $S$ is a simple semiring,

\item $S$ is a bounded distributive lattice.

\end{enumerate}

\begin{proof}

By Proposition \ref{simple}, the four statements $(1), (2), (3)$, and $(4)$ are equivalent. Obviously $(5)$ implies $(4)$. Now we prove that $(4)$ implies $(5)$.

$(4) \Rightarrow (5)$: In order to prove that $S$ is a bounded distributive lattice, we only need to prove the distributivity of addition on multiplication and the two absorption laws:

Distributivity: $(s+t)(s+u) = s^2+su+st+tu = s+su+st+tu = s(1+u)+st+tu = s+st+tu = s(1+t)+tu = s+tu$.

Absorption 1: $s+st = s(1+t) = s$.

Absorption 2: $s (s+t) = s^2 + st = s+st = s$.
\end{proof}

\end{proposition}

\begin{corollary}

\label{pbdl}

Let $S$ be a pseudocomplemented semiring. Then the following statements are equivalent:

\begin{enumerate}

\item $S$ is multiplicatively idempotent and $1\in \Stone(S)$,

\item $S$ is a bounded distributive lattice.

\end{enumerate}

\end{corollary}

Later in Theorem \ref{minimalpbdl}, we will discuss minimal primes of bounded distributive lattices with pseudocomplementation.

\begin{proposition}

\label{stone2}

Let $S$ be a pseudocomplemented semiring. Then the following statements are equivalent:

\begin{enumerate}

\item $S$ is a Stone semiring.

\item $(st)^* = s^* + t^* $ for any $s,t \in S$ and $1\in \Stone(S)$.

\end{enumerate}

\begin{proof}

$(1) \Rightarrow (2)$: It is clear that $st (s^* + t^*) = 0$. Now we show that $s^* + t^*$ is the largest element that annihilates $st$. Let $stx=0$ for some $x\in S$. Then by Proposition \ref{pseudo1}, $s^*{^*} tx =0$, which implies that $x s^*{^*} \leq t^*$. On the other hand, $xs^* \leq s^*$. Therefore we have the following: $x = x \cdot 1 = x(s^* + s^*{^*}) = xs^* + xs^*{^*} \leq s^* + t^*$.

$(2) \Rightarrow (1)$: $s^* + s^*{^*} = (ss^*)^* = 0^* = 1$.
\end{proof}

\end{proposition}

\begin{theorem}

\label{skel1}

Let $S$ be a Stone semiring. Then $(\Skel(S), \vee, \wedge)$ is a bounded complemented lattice, where $s \vee t = (s^* t^*)^*$ and $s \wedge t = s t$.

\begin{proof}

Let $S$ be a Stone semiring. Since $1\in \Stone(S)$, we have that $0^* =1$, which means that $1\in \Skel(S)$ and for any $s\in \Skel(S)$, $0 \leq s \leq 1$ (Proposition \ref{pseudo1}).

Now let $s,t \in \Skel(S)$. Obviously $s^*{^*} = s$. It is clear that $st \leq s$ and so, $(st)^*{^*} \leq s^*{^*} = s$. In a similar way, $(st)^*{^*} \leq t$. Note that every element of $S$ is pseudocomplemented and therefore $(st)^*{^*}\in \Skel(S)$. But since $\Stone(S) = S$, by Proposition \ref{stone1}, each element of $\Skel(S)$ is multiplicatively idempotent. From this we get that $(st)^*{^*} = (st)^*{^*} (st)^*{^*} \leq st$. But clearly, $st \leq (st)^*{^*}$. This means that $(st)^*{^*} = st$, i.e., $st \in \Skel(S)$. Finally let $x\in \Skel(S)$ such that $x\leq s,t$. It is now clear that $x = xx \leq st$. This means that $st = \inf_{\Skel(S)} \{s,t\}$.

On the other hand, $s^* t^* \leq s^*$. This means that $s = s^*{^*} \leq  (s^* t^*)^*$. In a similar way, $t \leq  (s^* t^*)^*$. Now let $x\in \Skel(S)$ such that $s,t \leq x$. Then $x^* \leq s^*, t^*$. So we have that $x^* = x^* x^* \leq s^* t^*$. This implies that $ (s^* t^*)^* \leq x$. This means that $\sup_{\Skel(S)} \{s,t\} = (s^* t^*)^*$.

From all we said we get that $(\Skel(S), \vee, \wedge)$ is a bounded lattice, where $s \vee t = (s^* t^*)^*$ and $s \wedge t = s t$. Also note that $s \wedge s^* = s s^* = 0$ and $s \vee s^* = (s^* s^*{^*})^* = 0^* = 1$. This already means that the lattice $(\Skel(S), \vee, \wedge)$ is complemented and the proof is complete.
\end{proof}

\end{theorem}

Let us recall that a semiring $S$ is complemented if for any $s\in S$, there exists $s^* \in S$ such that $ss^* = 0$ and $s + s^* = 1$.

\begin{corollary}

\label{Skelboolean}

Let $S$ be a pseudocomplemented semiring. If $S$ is a Stone semiring, then $(\Skel(S), +, \cdot, 0 , 1, ^*)$ is a boolean algebra. Conversely, if $(\Skel(S), +, \cdot, 0 , 1, ^{\prime})$ is a boolean algebra, then $S$ is a Stone semiring.

\begin{proof}

$\Rightarrow$: Let $S$ be a Stone semiring. Since $1^* = 0$, $0^* = 1$, and $0,1, 0^*, 1^* \in \pcomp(S)$, it is clear that $0,1 \in \Skel(S)$. Let $s,t \in \Skel(S)$. Then by Proposition \ref{stone2} and Theorem \ref{skel1}, $s \vee t = (s^* t^*)^* = s^*{^*} + t^*{^*} = s+t$. This means that $s+t \in \Skel(S)$. But as we have seen in the proof of Theorem \ref{skel1}, $st \in \Skel(S)$. All these observations assert that $\Skel(S)$ is a subsemiring of $S$. Now let $s\in \Skel(S)$. It is clear that $ss^* = 0$ and $s + s^* = s^*{^*} + s^*= 1$. This shows that any element of $\Skel(S)$ is complemented. On the other hand, since $s+s = s\vee s = s$ and $ss = s\wedge s = s$, $S$ is an idempotent semiring. Finally, since $0 \leq s \leq 1$ for any $s\in \Skel(S)$, $s+1 = 1$ for any $s\in \Skel(S)$, which means that $\Skel(S)$ is a simple semiring and by Proposition \ref{bdl}, $(\Skel(S), +, \cdot, 0 , 1, ^*)$ is a boolean algebra.

$\Leftarrow$: Let $(\Skel(S), +, \cdot, 0 , 1, ^{\prime})$ be a boolean algebra. Also let $s$ be an arbitrary element of $S$. Clearly $s^* \in \Skel(S)$ and $s^* (s^*)^{\prime} = 0$. This implies that $(s^*)^{\prime} \leq s^*{^*}$. From this we get that $1= s^* + (s^*)^{\prime} \leq s^* + s^*{^*}$. Clearly this implies that $s^* + s^*{^*} = 1$. This means that $S$ is a Stone semiring and the proof is complete.
\end{proof}

\end{corollary}

\begin{proposition}

\label{mi}

If a semiring $S$ is multiplicatively idempotent and positive, then we have the following:

\begin{enumerate}

\item If $s,t \in \pcomp(S)$, then $s+t \in \pcomp(S)$ and $(s+t)^* = s^* t^*$.

\item If $s, s^* \in \pcomp(S)$, then $(s+s^*)^* = 0$.

\item If $s,t,s^*, t^*, s^*{^*}, t^*{^*} \in \pcomp(S)$, then $(s^*{^*} + t^*{^*})^* = (s+t)^*$.

\item If $s,t \in \Skel(S)$, then $st \in \Skel(S)$.

\end{enumerate}

\begin{proof}
(1): Let $s,t \in \pcomp(S)$. It is clear that $s^* t^* (s+t)= 0$. Now let $u(s+t) = 0$. This point that $0$ is the least element of $S$, implies that $us = ut = 0$. This means that $u \leq s^*, t^*$. Now we observe that $u = u^2 \leq u t^* \leq s^* t^*$, which means that $s^* t^*$ is the largest element of $S$ that annihilates $s+t$. This means that $s+t \in \pcomp(S)$ and $(s+t)^* = s^* t^*$.

(2): Since $s, s^* \in \pcomp(S)$, by (1), we have $(s + s^*)^* = s^* s^*{^*} = 0$.

(3): Since $s^*{^*}, t^*{^*} \in \pcomp(S)$, $s^*{^*} + t^*{^*} \in \pcomp(S)$ and $(s^*{^*} + t^*{^*})^* = s^*{^*}{^*} t^*{^*}{^*} = s^* t^* = (s+t)^* $.

(4): Let $s,t \in \Skel(S)$. Then there are $u,v \in \pcomp(S)$ such that $s = u^*$ and $t= v^*$. Now by (a), we have that $st = u^* v^* = (u+v)^*$ and $u+v \in \pcomp(S)$. This implies that $st \in \Skel(S)$.
\end{proof}

\end{proposition}

We finalize this section by defining dense elements of a semiring:

\begin{definition}
Let $S$ be a positive semiring. We define $s\in \pcomp(S)$ to be a dense element of $S$ if $s^* = 0$. We denote the set of all dense elements of a semiring $S$ by $\Dns(S)$.
\end{definition}

\begin{proposition}

\label{dense1}

Let $S$ be a positive semiring. Then the following statements hold:

\begin{enumerate}

\item $1 \in \Dns(S)$.

\item If $s,t \in pcomp(S)$, $s\leq t$, and $s\in \Dns(S)$, then $t \in \Dns(S)$.

\end{enumerate}

If in addition the semiring $S$ is multiplicatively idempotent, we have the following:

\begin{enumerate}

\item[(a)] If $s \in \pcomp(S)$ and $t \in \Dns(S)$, then $s+t \in \Dns(S)$.

\item[(b)] If $s, s^* \in \pcomp(S)$, then $s+s^* \in \Dns(S)$.

\item[(c)] If $s,t,s^*, t^*, s^*{^*}, t^*{^*} \in \pcomp(S)$, then $s^*{^*} + t^*{^*} \in \Dns(S)$ if and only if $s+t \in \Dns(S)$.

\end{enumerate}

\begin{proof}
The proof of (1) and (2) is straightforward and the statements (a), (b), and (c) are obtained from Proposition \ref{mi}.
\end{proof}

\end{proposition}

\section{Minimal Prime Ideals of Semirings}

Let us recall that a nonempty subset $I$ of a semiring is called an ideal if $a+b \in I$ and $sa \in I$ for all $a,b \in I$ and $s\in S$. An ideal $P \neq S$ is defined to be a prime ideal of $S$, if $ab\in P$ implies either $a\in P$ or $b\in P$. A prime ideal $P$ of a semiring $S$ is called to be a minimal prime ideal of $S$, if $I \subseteq P$ implies either $I = (0)$ or $I = P$. For more on ideals of a semiring, one can refer to chapter 6 and chapter 7 of the book \cite{Golan1999}. A nonempty subset $W$ of a semiring $S$ is said to be a multiplicatively closed set (for short an MC-set) if $1\in W$ and for all $w_1,w_2 \in W$, we have $w_1 w_2 \in W$. In other words, $W$ is an MC-set if and only if it is a submonoid of $(S,\cdot)$. It is clear that an ideal $P$ of $S$ is a prime ideal of $S$ if and only if $S-P$ is an MC-set. The following theorem is semiring version of a theorem in commutative algebra due to German mathematician Wolfgang Krull (1899-1971):

\begin{theorem}

\label{maxisprime1}

The maximal elements of the set of all ideals disjoint from an MC-set of a semiring are prime ideals.

\begin{proof}
The proof is just a mimic of the proof of \cite[Theorem 1, p. 1]{Kaplansky1970} and therefore omitted.
\end{proof}

\end{theorem}

The ring version of the following theorem is also credited to Krull:

\begin{theorem}

\label{krullintersectionprimes}

Let $S$ be a semiring and $I$ an ideal of $S$. Then $\sqrt{I} = \bigcap_{P\in V(I)} P$, where $V(I) = \{P\in \Spec(S) : P \supseteq I\}$.

\begin{proof}
It is straightforward that $\sqrt{I} \subseteq \bigcap_{P\in V(I)} P$. Now let $s\notin \sqrt{I}$. It is clear that $W_s = \{s^n: n\geq 0\}$ is an MC-set of $S$ disjoint from $\sqrt{I}$. So there exists a prime ideal containing $I$ and not containing $s$.
\end{proof}

\end{theorem}

 Now if we consider $V(I) = \{P\in \Spec(S): P \supseteq I\}$, partially ordered by containment relation, then by Zorn's Lemma it has a $\supseteq$-maximal element, which is a $\subseteq$-minimal element. Those prime ideals, which are $\subseteq$-minimal elements of $V(I)$ are called minimal primes of $I$ and one may collect them in a set denoted by $\Min(I)$. Usually $\Min(0)$ is denoted by $\Min(S)$ and the elements of $\Min(S)$ are called minimal primes of $S$. A semiring $E$ is called to be entire, if $ab=0$ implies that either $a=0$ or $b=0$. It is, by definition, clear that if $E$ is an entire semiring, then $(0)$ is the only minimal prime of $E$. In this case, minimal elements of the set of nonzero prime ideals of the semiring $E$ may play an important role and are usually called height 1 prime ideals of $E$.

\begin{corollary}

\label{krullintersectionprimes2}

Let $S$ be a semiring and $I$ an ideal of $S$. Then $\sqrt{I} = \bigcap_{P\in \Min(I)} P$, where by $\Min(I)$ we mean the set of all minimal primes of $I$.

\end{corollary}

\begin{remark}

An element $s\in S$ is said to be nilpotent if $s^n = 0$ for some $n\in \mathbb N$. The set of all nilpotent elements of the semiring $S$ is called the \emph{lower nil radical} of $S$ and is denoted by $\Nil (S)$. It is clear that $\Nil (S) = \sqrt{(0)}$. We call a semiring $S$ to be nilpotent-free if $\sqrt{(0)} = (0)$. In other words, a semiring $S$ is nilpotent-free if the only nilpotent element of $S$ is the element $0$. While this condition in ring theory is known as ``reduced" (Cf. \cite[p. 3]{Matsumura1989}), we prefer not to use this term for this concept, since it has been reserved for another concept in semiring theory (Cf. \cite[Example 8.8]{Golan1999}).

\end{remark}

The following theorem and its corollaries are the semiring version of Theorem 2.1 in \cite{Huckaba1988} and its corollaries. For more on minimal prime ideals of reduced rings, one can also see \cite{Matlis1983}.

\begin{theorem}

\label{minimalprimehuckaba}

Let $P \supseteq I$ be ideals of a semiring $S$, where $P$ is prime. Then the following statements are equivalent:

\begin{enumerate}

\item $P$ is a minimal prime ideal of $I$.

\item $S-P$ is an MC-set maximal with respect to missing $I$.

\item For each $x\in P$, there is a $y\notin P$ and a nonnegative integer $i$ such that $yx^i \in I$.

\end{enumerate}

\begin{proof}

$(1) \Rightarrow (2)$: Expand $S-P$ to an MC-set, say $W$, which is maximal with respect to missing $I$. Let $Q$ be an ideal containing $I$ that is maximal with respect to being disjoint from $W$. By Theorem \ref{maxisprime1}, $Q$ is prime. But $Q \supseteq I$ is disjoint from $S-P$ and $P$ is a minimal prime of $I$, so $Q=P$ and therefore $W=S-P$.

$(2) \Rightarrow (3)$: Take a nonzero $x\in P$ and set $W := \{yx^i: y\in (S-P) \wedge i\in \mathbb N_0 \}$. Then $W$ is an MC-set that properly contains $S-P$. Therefore $W \cap I \neq \emptyset$, which means that there is a $y\notin P$ and a nonnegative integer $i$ such that $yx^i \in I$.

$(3) \Rightarrow (1)$: Suppose $I \subset Q \subseteq P$, where $Q$ is a prime ideal of $S$. If there is some $x\in P-Q$, then there is a $y\notin P$ and a positive integer $i$ such that $yx^i \in I \subset Q$, a contradiction. So $Q = P$ and this finishes the proof.
\end{proof}

\end{theorem}

\begin{corollary}

\label{minimalprimehuckaba2}

If $S$ is a nilpotent-free semiring and $P$ is a prime ideal of $S$, then $P$ is a minimal prime ideal of $S$ if and only if for each $x\in P$ there exists a $y\notin P$ such that $xy=0$.

\begin{proof}
By Theorem \ref{minimalprimehuckaba}, $\Leftarrow$ is obvious. Now let $P$ be a minimal prime ideal of $S$ and $x\in P$. If $x=0$, then there is nothing to prove. If $x \neq 0$, then there is a $y\notin P$ and a positive integer $i$ such that $xy^i =0$. This obviously implies that $(xy)^i =0$. But $S$ is nilpotent-free, so $xy=0$. Q.E.D.
\end{proof}

\end{corollary}

Let $S$ be a semiring and $H$ be a nonempty subset of $S$. The set of all annihilators of $H$, denoted by $\Ann(H):=\{s\in S : s\cdot H =(0)\}$ is an ideal of $S$. If $H=\{s\}$ is a singleton, then we write $\Ann(s)$ instead of $\Ann(\{s\})$. Particularly if $J = (s_1, \ldots, s_n)$ is a finitely generated ideal of $S$, then we simply write $\Ann(s_1, \ldots, s_n)$ instead of $\Ann(J)$. One can easily check that $\Ann(J) = \bigcap^n_{i=1} \Ann(s_i)$, whenever $J = (s_1, \ldots, s_n)$.

\begin{corollary}

\label{minimalprimehuckaba3}

Let $J$ be a finitely generated ideal of a nilpotent-free semiring $S$. Then $J$ is contained in a minimal prime ideal $P$ of $S$ if and only if $\Ann(J) \neq (0)$.

\begin{proof}

Let $S$ be a nilpotent-free semiring and $J = (s_1, \ldots,s_n)$ for some $s_1, \ldots,s_n \in S$.

$(\Rightarrow)$: If $J$ is contained in a minimal prime ideal $P$ of $S$, then for any $s_i$, there is a $t_i \in S-P$ such that $s_i t_i =0$. Take $t = t_1 \cdots t_n$. It is clear $t$ is a nonzero annihilator of $J$ and therefore $\Ann(J) \neq 0$.

$(\Leftarrow)$: Let $\Ann(J)\neq (0)$. Then by Theorem \ref{krullintersectionprimes2}, $\Ann(J)$ cannot be a subset of all minimal primes of $S$. So there is a minimal prime $P$ of $S$ such that $\Ann(J) \nsubseteq P$. Our claim is that $J \subseteq P$. In contrary let $J \nsubseteq P$. This implies that at least one of the generators of $J$, say $s_1$, is not an element of $P$. But by Corollary \ref{minimalprimehuckaba2}, for all $y\notin P$, we have $s_1 y \neq 0$. This means that $\Ann(s_1) \subseteq P$. But $\Ann(J) \subseteq \Ann(s_1)$, a contradiction. Consequently, $J \subseteq P$, the thing it was required to have shown.
\end{proof}

\end{corollary}

Let us recall that an element $s\in S$ is said to be a zero-divisor of the semiring $S$, if there is a nonzero element $t\in S$ such that $st =0$. The set of all zero-divisors of $S$ is denoted by $Z(S)$.

\begin{corollary}

\label{zero-divisorunionofminimalprimes}

If $S$ is a nilpotent-free semiring, then $Z(S) = \bigcup_{P\in \Min(S)} P$.

\begin{proof}
If $x\in P$ for some $P\in \Min(S)$, then by Corollary \ref{minimalprimehuckaba2}, $x$ has a nonzero annihilator and therefore $x\in Z(S)$. On the other hand, if $x\in Z(S)$. Then $\Ann(x) \neq (0)$ and therefore by Corollary \ref{minimalprimehuckaba3}, $(x) \subseteq P$ for some $P\in \Min(S)$ and the proof is complete.
\end{proof}

\end{corollary}

\begin{theorem}

\label{minimalpbdl}

Let $S$ be a mutiplicatively idempotent and pseudocomplemented semiring such that $1\in \Stone(S)$ and $P$ a prime ideal of $S$. Then the following statements are equivalent:

\begin{enumerate}

\item If $s\in P$, then $s^* \notin P$, for each $s\in S$,

\item If $s\in P$, then $s^*{^*} \in P$, for each $s\in S$,

\item $P \cap \Dns(S) = \emptyset$.

\item $P$ is a minimal prime ideal of $S$.

\end{enumerate}

\begin{proof}

$(1) \Rightarrow (2)$: Corollary \ref{minimalprimepseudo}.

$(2) \Rightarrow (3)$: Let $s\in P \cap \Dns(S)$, for some $s\in S$. So $s^*=0$ and therefore $s^*{^*} = 0^* = 1$, which implies that $1\in P$, a contradiction.

$(3) \Rightarrow (1)$: Let $s, s^* \in P$ for some $s\in S$. So $s+s^* \in P$. By Proposition \ref{dense1}, $ s+s^* \in \Dns(S)$, which means that $P \cap \Dns(S) \neq \emptyset$.

$(1) \Rightarrow (4)$: Corollary \ref{minimalprimepseudo}.

$(4) \Rightarrow (1)$: Let $s\in P$. Since $P$ is a minimal prime ideal of $S$. There is a $t \notin P$ such that $st=0$. This implies that $t \leq s^*$. Note that $S$ is a bounded distributive lattice by Corollary \ref{pbdl}. Now if $s^* \in P$, we have that $t = ts^* \in P$, a contradiction. So $s^* \notin P$ and the proof is complete.
\end{proof}

\end{theorem}

\section*{Acknowledgments}

This work is supported by University of Tehran. Our special thanks go to University of Tehran, College of Engineering and Department of Engineering Science for providing all the necessary facilities available to us for successfully conducting this research. The author of this work is grateful to Professor Henk Koppelaar for introducing the book \cite{RasiowaSikorski1963}.

\bibliographystyle{plain}

\end{document}